\documentclass{article}
\usepackage{amsthm,amsfonts, amsbsy, amssymb,amsmath,graphicx}
\usepackage{graphics}

\newtheorem{thm}{Theorem}
\newtheorem{lm}{Lemma}
\newtheorem{ex}{Example}
\newtheorem{crl}{Corollary}

\title{On Free Knots}

\author{Vassily Olegovich Manturov}

 \def\Z{{\mathbb Z}}

\begin{document}

\maketitle

\section{Statement of the Problem}

The aim of the present paper is to study {\em free knots}, a
dramatic simplification of the notion of virtual knots also
connected to finite type invariants, curves on surfaces and other
objects of low-dimensional topology. It turns out that the objects
of such sort (a factorization of arbitrary graphs modulo three
formal Reidemeister moves) can be studied by using some very much
simpler objects: in the paper we construct an invariant valued in
the set of equivalent classes of the same graphs modulo only {\em
the second Reidemeister move}. By using this invariant, one can
easily prove the non-triviality of some free knots; previously, the
{\em existence} of non-trivial free knots was not evident.

During the whole article by {\em graph} we mean a finite graph,
possibly, having loops and/or multiple edges, not necessarily
connected. Herewith, by a connected component of such graph we also
admit (besides graphs in the usual sense) free loops without
vertices. The number of connected components of a graph is thought
to be finite.

A $4$-valent graph $G$ (possibly, with free loops) is called {\em
framed} if for every vertex $V$ of $G$ the four emanating half-edges
are split into two pairs of (formally) opposite edges.

Denote by $G_{0}$ the one-component free loop.

By a {\em unicursal component} of a framed graph we mean the
equivalence class of its edges generated by the following elementary
equivalence relation: two edges are elementary equivalent if they
possess half-edges which attach a vertex from opposite sides.

Given a chord diagram \footnote{One may consider oriented or
non-oriented graphs (singular knots) with corresponding chord
diagrams being oriented or non-oriented. In this paper we restrict
ourselves for the non-oriented case, though the orientation can be
taken into account straightforwardly.} $C$. Then the corresponding
framed graph is constructed as follows. Edges of the graph $G(C)$
are assoicated with arcs of the chord diagrams, and vertices are
associated to chords of the chord diagram. Each vertex is incident
to four half-edges corresponding to arcs attaching two ends of a
chord. Arcs which are incident to the same chord end will correspond
to the (formally) opposite half-edges. The opposite procedure is
evident: having a framed graph with one unicursal component, one can
generate it by a certain chord diagram, to be denoted by $C(G)$.

Thus, there is a one-to-one correspondence between framed 4-graphs
with one unicursal component and chord diagrams.

For a given four-valent framed graph $G$ with one unicursal compnent
with every two vertices $v_{1},v_{2}$ we associate an element
$\langle v_{1},v_{2}\rangle$ from ${\bf Z}_{2}=\{0,1\}$ which is
equal to one if the corresponding chords of $C(G)$ are unlinked
(recall that two chords $a,b$ of a chord diagram $C$ are  {\em
linked} if after removeing the two ends of $a$ from the circle of
the chord diagram, the ends of the chord $b$ lie in different
connected components.

Later we assume that all framed graphs have one unicursal component,
unless otherwise specified.

A simplest example of a framed graph is a knot projection. This
gives a {\em planar framed graph}. More generally, one can consider
a virtual knot diagram where classical crossings play the role of
vertices, and virtual crossings are just intersection points of
images of different edges.

It is known that classical and virtual knots are encoded by {\em
Gauss diagrams}, see \cite{GPV}. The notion of framed graph is a
serious simplification of the notion of Gauss diagram. First, at
classical vertices we don't indicated which branch forms an overpass
and which one forms an underpass; besides, we do not indicate the
writhe number of the vertex; thus, we forget both arrows and signs
of the Gauss diagram.

Therefore, a framed graph is a very strong simplification of a knot
diagram.

Our goal is to consider framed 4-graphs modulo some moves analogous
to Reidemeister moves and study the obtained graph equivalence
classes.

By a {\em free knot} we mean an equievalence class of framed
$4$-valent graphs with one unicursal component modulo the following
transformations. For each transformation we assume that only a fixed
fragment of the graph is being operated on (this fragment is to be
depicted) or some corresponding fragments of the chord diagram. The
remaining part of the graph or chord diagram are not shown in the
picture; the pieces of the chord diagram not containing chords
participating in this transformation, are depicted by punctured
arcs. The parts of the graph are always shown in a way such that the
formal framing (opposite edge relation) in each vertex coincides
with the natural opposite edge relation taken from  ${\bf R}^{2}$.

\begin{figure}
\centering\includegraphics[width=200pt]{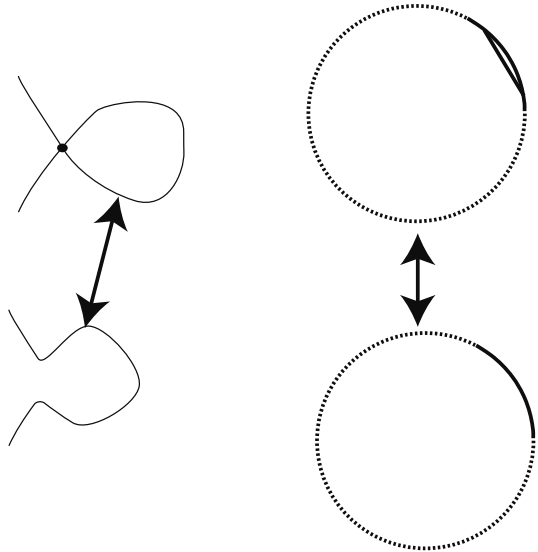}
\caption{Addition/removal of a loop on a graph and on a chord
diagram} \label{1r}
\end{figure}

The first Reidemeister move is an addition/removal of a loop, see
Fig.\ref{1r}

\begin{figure}
\centering\includegraphics[width=300pt]{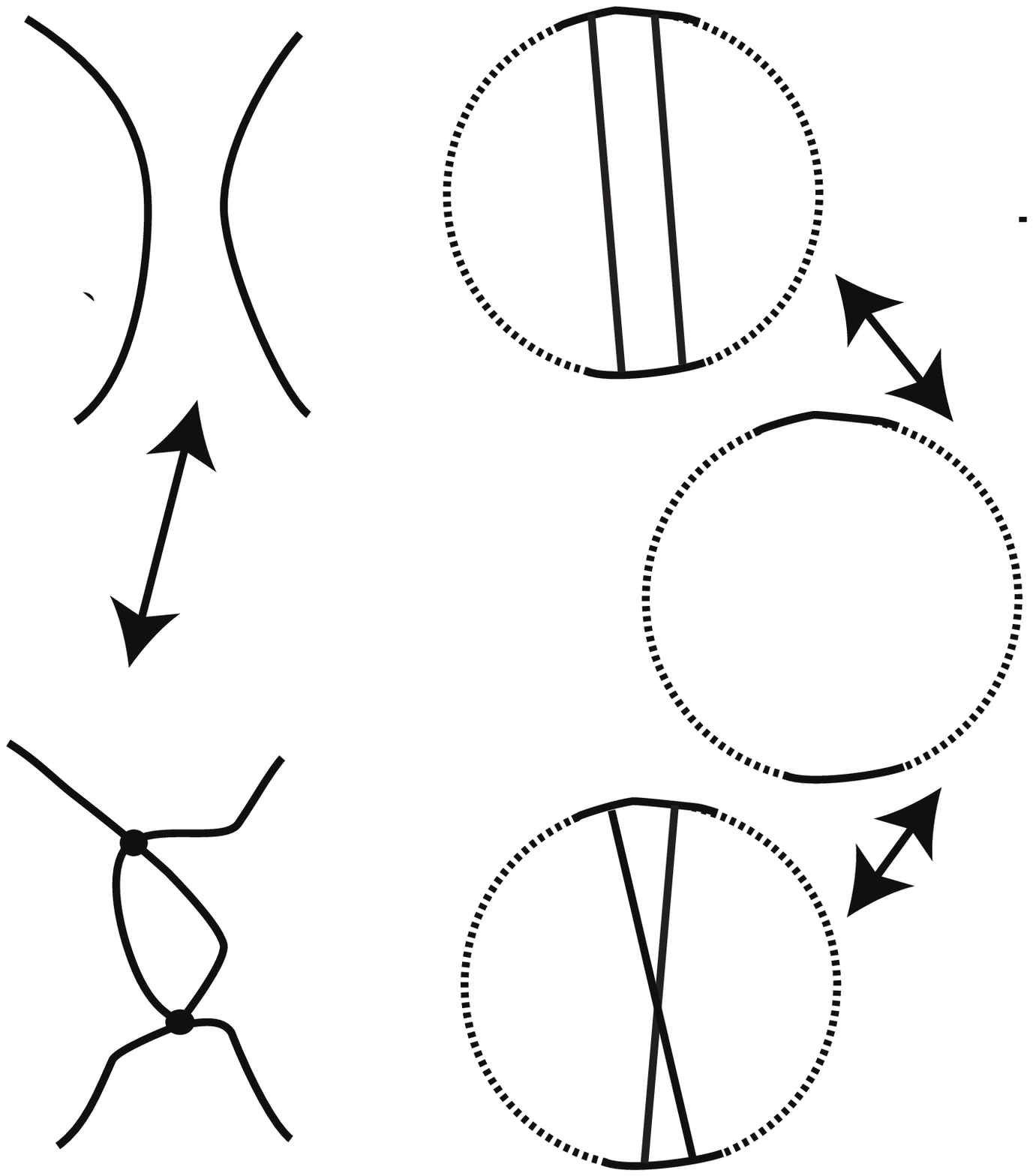} \caption{The second
Reidemeister move and two chord diagram versions of it} \label{2r}
\end{figure}

The second Reidemeister move adds/removes a bigon formed by a pair
of edges which are adjacent in two edges, see Fig. \ref{2r}.

Note that the second Reidemeister move adding two vertices does not
impose any conditions on the edges it is applied to.

\begin{figure}
\centering\includegraphics[width=300pt]{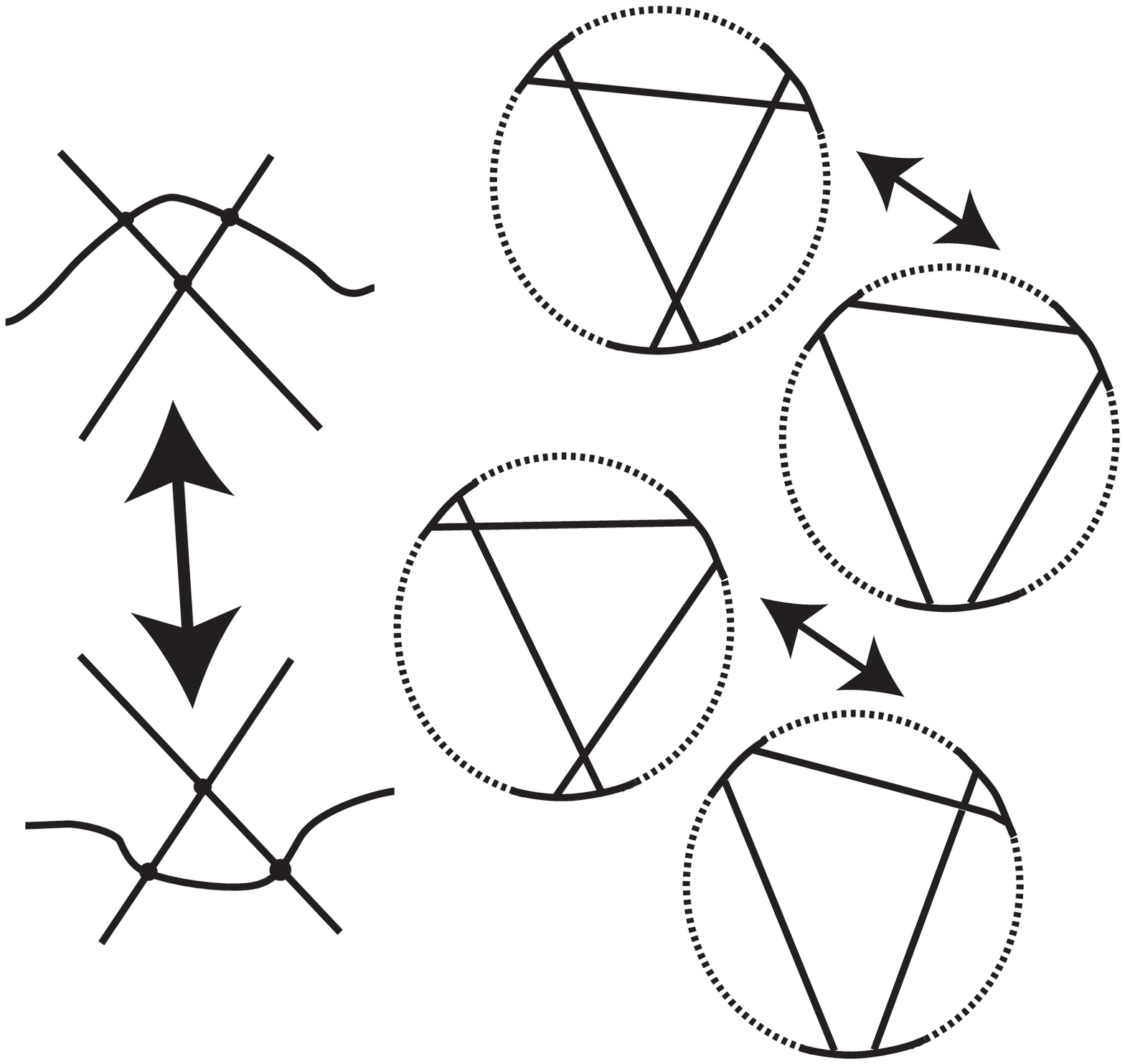} \caption{The third
Reidemeister move and its chord diagram versions} \label{3r}
\end{figure}

The third Reidemeister move is shown in Fig.\ref{3r}.

Note that each of these three moves applied to a framed graph,
preserves the number of unicursal components of the graph. Thus,
applying these moves to graphs with a unique unicursal cycle, we get
to graphs with a unique unicursal cycle.

A {\em free knot} is an equivalence class of framed graphs (with a
unique unicursal cycle) modulo the transformations listed above.

{\em What are free knots?} If we consider a planar framed graph
(which originates from a classical knot) then this planar graph can
be easily reduced to $G_0$. It can be easily shown that a
one-component framed graph embeddable in torus (with framing
preserved) is also reducible to $G_0$.

Free knots are closely connected to  {\em flat virtual knots},
i.\,e., with equivalence classes of virtual knots modulo
transformation changing over/undercrossing structure. The latter are
equivalence classes of immersed curves in orientable $2$-surfaces
modulo homotopy and stabilization.

Nevertheless, the equivalence of free knots is even stronger than
the equivalence of flat virtual knots: our free knots do not require
any surface. Every time one applies a Reidemeister move to a regular
4-graph, one embeds this graph into a 2-surface arbitrarily (with
framing preserved), apply this Reidemeister move inside the surface
and then forget the surface again.

\begin{ex}
Consider the flat virtual Kishino knot, see Fig. \ref{Kishino}

\begin{figure}
\centering\includegraphics[width=200pt]{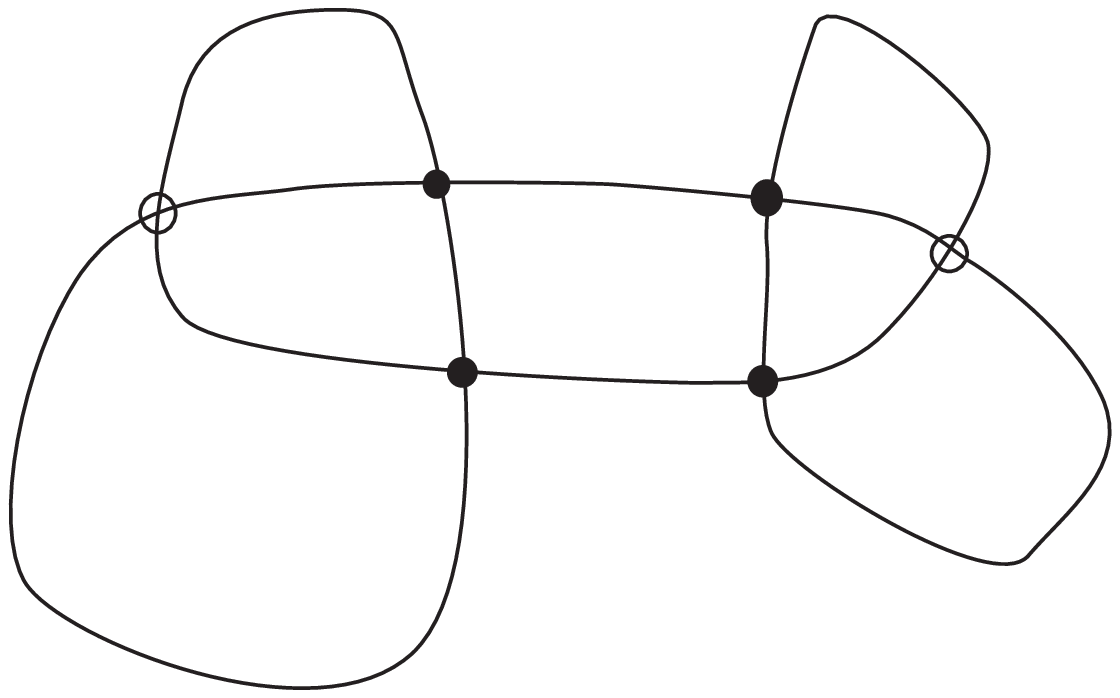} \caption{The
Flat Kishino knot} \label{Kishino}
\end{figure}
\end{ex}

It is known that this not is not trivial as a flat virtual knot: the
corresponding representative lies on the sphere with two handles and
it is minimal (this representative splits the sphere into a
quadrilateral and an octagon).

Nevertheless, the corresponding framed graph considered by itself
has a bigon formed by two edges which are adjacent in two vertices.
Thus, the free knot represented by the flat Kishino knot is trivial.

The exact statement connecting virtual knots and free knots sounds
as follows:

\begin{lm}
A free knot is an equivalence class of virtual knots modulo two
transformations: crossing switches and virtualizations.

A virtualization is a local transformation shown in Fig.
\ref{virtua}.

\begin{figure}
\centering\includegraphics[width=200pt]{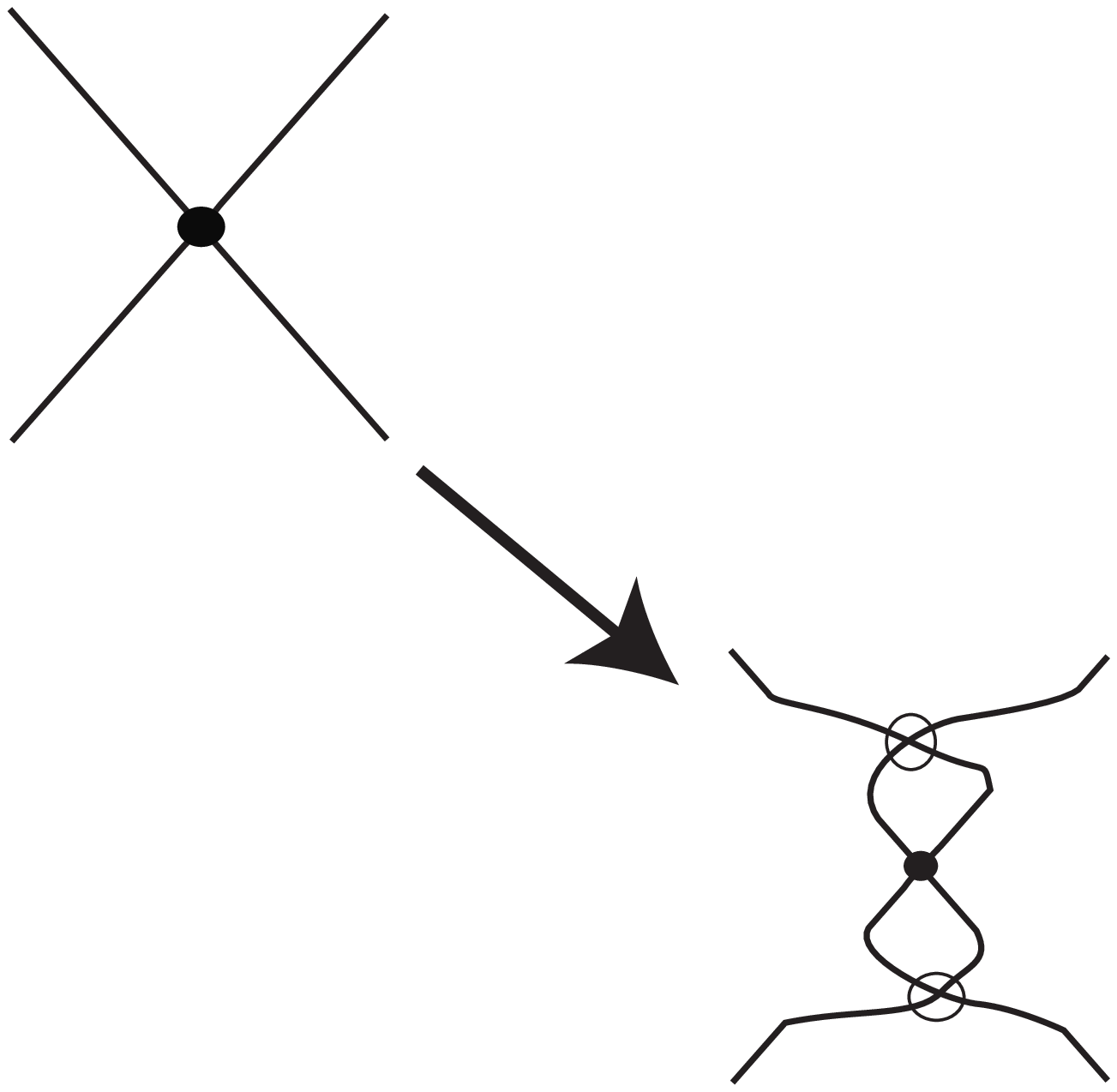}
\caption{Virtualization} \label{virtua}
\end{figure}
\end{lm}

One may think of a virtualization as way of changing the immersion
of a 4-valent framed graph in plane.

Note that in the case of {\em free links},i.\,e., framed 4-graphs
with more than one unicursal components, it is much easier to find a
non-trivial example.

Consider the framed graph $G_2$ with one vertex $x$ and two edges
$a,b$, each connecting  $x$ to $x$ in such a way that the edge $a$
is opposite to $a$, and the edge $b$ is opposite to $b$.

This {\em free link} is not equivalent to the trivial one because of
the following simple {\em invariant} of two-component free links.
Consider a 2-component free link and calculate the {\em parity} of
the number of crossings formed by two components. This parity is
obviously invariant under Reidemeister moves.

Thus,  $G_2$ is not equivalent to the unlink (formed by two disjoint
free loops).

In the case of knots we shall also seek some {\em oddness} which can
not be got rid of by using Reidemeister moves.

Given a framed $4$-graph $G$ (with one unicursal component), and let
$V=\{v_{1},\dots, v_{n}\}$ be the set of vertices of $G$. Consider
the linear ${\Z}_{2}$-space $L(G)$ generated by vectors
$v_{1},\dots,v_{n}$. We shall identify subsets of the set of
vertices of $G$ with vectors from $L(G)$ and we shall consider sums
of elements from $L(G)$ as boolean sums in $2^{V}$.

A vertex of $G$ is called {\em even}, if the number of vertices
incident to it is even, and {\em odd} otherwise.

For a vertex $a$, we denote by $E_{a}$  the sum of vertices incident
to $a$ (by definition, we assume that the vertex is not incident to
itself).

We emphasize some evident properties which hold under Reidemeister
moves.

The first Reidemeister move corresponds to an addition/removal of an
even vertex, and the pairwise incidence relation of the remanining
vertices does not change.

The second Reidemeister move adds (removes) two vertices $a,b$ of
the same parity, herewith  $E_{a}+E_{b}$ is either $0$ or coincides
with  $a+b$; after applying the Reidemeister move the parity of the
remaining vertices does not change. Neither does the pairwise
incidence of the remaining vertices.

Note that for every chord diagram, the chords are naturally split
into {\em equivalence classes}: $a\sim b$ if and only if
$E_{a}+E_{b}=0$ or $E_{a}+E_{b}=a+b$. We call a set of equivalence
classes of chords a {\em bunch}; in every bunch, all chords are
either pairwise linked or pairwise unlinked. Every second
Reidemeister move adds/deletes two chords from the same bunch. Note
that this does not affect the bunches (equivalences of the remaining
chords).

Note also that the intersection index of two chords from different
bunches $\alpha,\beta$ does not depend on the particular choice of
these chords, thus one can write $\langle \alpha,\beta\rangle$.

When performing the third Reidemeister move, we have three vertices
(chords) $a,b,c$, for which  $E_{a}+E_{b}+E_{c}\subset \{a,b,c\},
|E_{a}+E_{b}+E_{c}|=0$ or $2$. After performing the move, we get
instead of $a,b,c$ three vertices $a',b',c'$ with pairwise switched
incidences (w.r.t. $a,b,c$) the remaining incidences are unchanged:
for $f,g\notin \{a,b,c\}: \langle f,g\rangle$ remains unchanged an
$\langle f,a\rangle=\langle f,a'\rangle;\langle f,b\rangle=\langle
f,b'\rangle;\langle f,c\rangle=\langle f, c'\rangle$. Therefore it
follows that, e.\,g., the number of odd vertices amongst $a,b,c$ is
{\bf even} (is equal to zero or two).

Herewith the parity of $a$ coincides with that of $a'$, the parity
of $b$ coincides with that of $b'$, and the parity of  $c$ coincides
with that of $c'$

Thus, to decrease the number of vertices of the diagram by using
Reidemeister move, it is necessary to have either an even vertex (to
be able to apply the first Reidemeister move) or a couple of
vertices having the same incidence with any of the remaining
vertices. The third Reidemeister move can be applied if for some
vertices $a,b,c$ we have $E_{a}+E_{b}+E_{c}\subset \{a,b,c\}$.

Thus, a natural class of graphs arises, where neither a
simplification nor a third Reidemeister move can be applied in turn.

For example, so are the graphs with all vertices odd and such that
for every two distinct vertices $a,b$ there exists a vertex
$c\notin\{a,b\}$ such that $\langle a,c\rangle\neq \langle
b,c\rangle$.

We call such graphs {\em irreducibly odd}.

\begin{figure}
\centering\includegraphics[width=200pt]{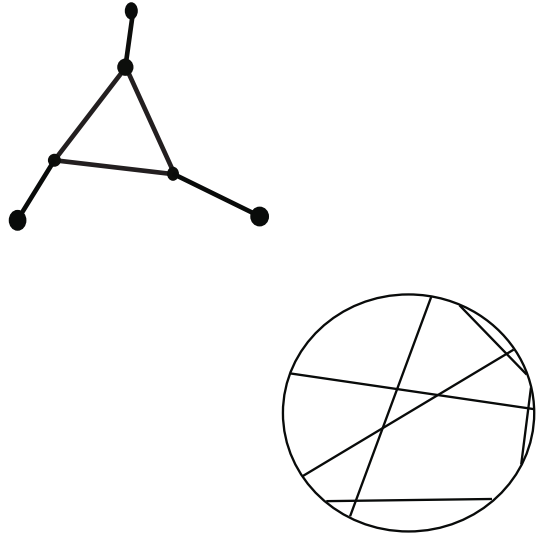} \caption{An
irreducibly odd graph and its chord diagram} \label{irred}
\end{figure}

The simplest example of an irreducibly odd graph is depicted in Fig.
\ref{irred}.

Assume an irreducibly odd graph $G$ generates a free knot $K$. By
definition it is impossible to simplify $G$ in one turn. It turns
out that the representative $G$ of the knot $K$ is indeed {\em
minimal}: any other representative of $K$ has the number of vertices
at least as many as those of $G$.

{\bf The main idea behind the proof is that in certain cases the
equivalence recognition problem with respect to all three
Reidemeister moves can be reduced to the equivalence recognition
problem with respect to only Reidemeister 2 move, which is
significantly easier}.

Now, we are going to construct an invariant map of free knots valued
in some objects considered modulo only the second Reidemeister
moves.

\section{The Main Theorem}

\newcommand{\ZG}{\Z_{2}{\mathfrak{G}}}
Let ${\mathfrak{G}}$ be the set of all equivalence classes of framed
 graphs with one unicursal component modulo second Reidemeister moves.
Consider the linear space $\ZG$.

Let $G$ be a framed graph, let  $v$ be a vertex of $G$ with four
incident half-edges $a,b,c,d$, s.t.  $a$ is opposite to $c$ and $b$
is opposite to $d$ at $v$.

\begin{figure}
\centering\includegraphics[width=200pt]{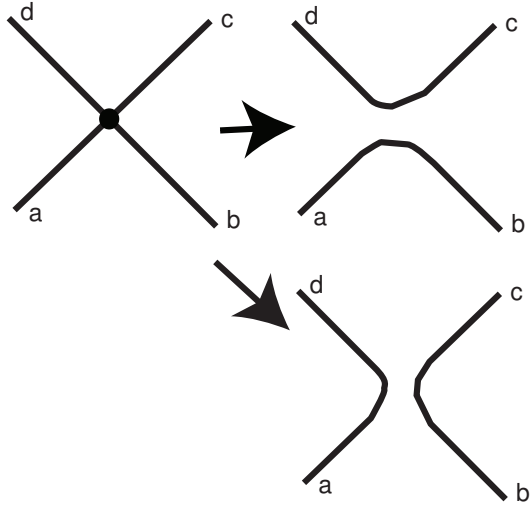} \caption{Two
smoothings of a vertex of for a framed graph} \label{smooth}
\end{figure}

By {\em smoothing} of $G$ at $v$ we mean any of the two framed
$4$-graphs obtained by removing $v$ and repasting the edges as
$a-b$, $c-d$ or as $a-d,b-c$, see Fig. \ref{smooth}.

Herewith, the rest of the graph (together with all framings at
vertices except $v$) remains unchanged.

We may then consider further smoothings of $G$ at {\em several}
vertices.

Consider the following sum

\begin{equation}
[G]=\sum_{s\;even.,1\; comp} G_{s},
\end{equation}
which is taken over all smoothings in all {\em even} vertices, and
only those summands are taken into account where $G_{s}$ has one
unicursal component.

Thus, if  $G$ has $k$ even vertices, then $[G]$ will contain at most
$2^{k}$ summands, and if all vertices of $G$ are odd, then we shall
have exactly one summand, the graph $G$ itself.

Consider  $[G]$ as an element of $\ZG$. In this case it is evident,
for instance, that if all vertices of $G$ are even then $[G]=[G_0]$:
by construction, all summands in the definition of $[G]$ are equal
to $[G_0]$, it can be easily checked that the number of such
summands is odd.

Now, we are ready to formulate the main theorem:

\begin{thm}
If $G$ and $G'$ represent the same free knot then in $\ZG$ the
following equality holds: $[G]=[G']$.\label{mainthm}
\end{thm}

Theorem \ref{mainthm} yields the following
\begin{crl}
Let  $G$ be an irreducibly odd framed 4-graph with one unicursal
component. Then any representative $G'$ of the free knot
$K_{G}$,generated by $G$, has a smoothing  $\tilde G$ having the
same number of vertices as $G$. In particular, $G$ is a minimal
representative of the free knot $K_{G}$ with respect to the number
of vertices.\label{sld}
\end{crl}

\subsection{The set ${\ZG}$}.

Having a framed $4$-graph, one can consider it as an element of
$\ZG$. It is natural to try simplifying it: we call a graph in $\ZG$
{\em irreducible} if no decreasing second Reidemeister move can be
applied to it. We call a graph in $\ZG$ {\em irreducible} if it has
no free loops and no decreasing second Reidemeister move can be
applied to it.

The following theorem is trivial

\begin{thm}
Every $4$-valent framed graph $G$ with one unicursal component
considered as an element of $\ZG$ has a unique irreducible
representative, which can be obtained from $G$ by consequtive
application of second decreasing Reidemeister moves.

\end{thm}

This allows to recognize elements $\ZG$  easily, which makes the
invariants constructed in the previous subsection digestable.

In particular, the minimality of a framed $4$-graph in $\ZG$  is
easily detectable: one should just check all pairs of vertices and
see whether any of them can be cancelled by a second Reidemeister
move (or in $\ZG$ one should also look for free loops)
.
\begin{proof}[Proof of the Corollary]
By definition of $[G]$ we have $[G]=G$. Thus if $G'$ generates the
same free knot as $G$ we have $[G']=G$ in $\ZG$.

Consequently, the sum representing $[G']$ in $\mathfrak{G}$ contains
at least one summand which is $a$-equivalent to $G$. Thus $G'$ has
at least as many vertices as $G$ does.

Moreover, the corresponding smoothing of $G'$ is a diagram, which is
 $a$-equivalent to $G$. One can show that under some (quite natural)
``rigidity'' condition this will yield that one of smoothings of
$G'$ coincides with $G$.
\end{proof}

\begin{proof}[Proof of the Theorem]
Let us check the invariance $[G]\in \ZG$ under the three
Reidemeister moves.

Let  $G'$ differ from  $G$ by a first Reidemeister move, so that
$G'$ has one vertex more than $G$. By definition this vertex is
even, and when calculating $[G']$ this vertex has to be smoothed in
order to get one unicursal curve in total.

Thus, we have to take only one of two smoothings of the given
vertex, see Fig. \ref{razved}

\begin{figure}
\centering\includegraphics[width=200pt]{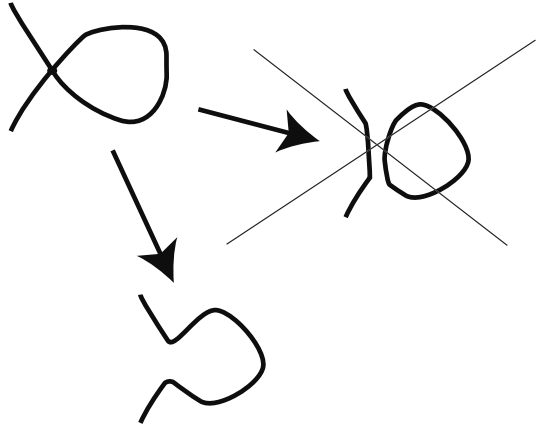} \caption{The two
smoothings --- good and bad --- of a loop} \label{razved}
\end{figure}

Thus there is a natural equivalence between smoothings of $G$ having
one unicursal component, and smoothings of $G'$ with one unicursal
component. Moreover, this equivalence yields a termwise identity
between $[G]$ and $[G']$.

Now, let  $G'$ be obtained from $G$ by a second Reidemeister move
adding two vertices.

These two vertices are either both even or both odd.

If both added vertices are odd, then the set of smoothings of $G$ is
in one-to-one correspondence with that of $G'$ and the corresponding
summands for $[G]$ and for $[G']$ differ from each other by a second
Reidemeister move.

If both vertices are odd then one has to consider different
smoothings of these vertices shown in. \ref{razved2}.

\begin{figure}
\centering\includegraphics[width=300pt]{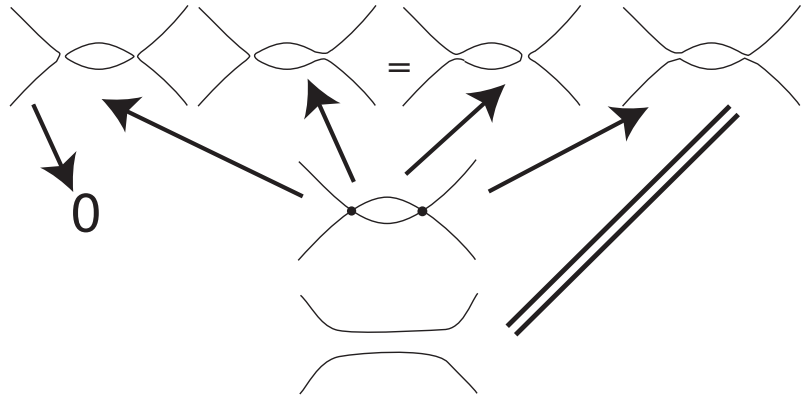}
\caption{Smoothings of two even vertices} \label{razved2}
\end{figure}

The smoothings shown in the upper-left Fig. \ref{razved2}, yield
more than one unicursal component (there is a free loop), so they do
not count in $[G']$.

The second-type and third-type smoothings (the second and the third
pictures in the top row of \ref{razved2}) give the same impact to
 $\ZG$, thus, they reduce in
$[G']$. Finally, the smoothings corresponding to the upper-right
Fig. \ref{razved2} are in one-to-one correspondence with smoothings
of $G$, thus we have a term-wise equality of terms $[G]$ and those
terms of $[G']$, which are not cancelled by comparing the two middle
pictures.

If $G$ and $G'$ differ by a third Reidemeister move, then the
following two cases are possible: either all vertices taking part in
the third Reidemeister move are even, or two of them are odd and one
is even.

If all the three vertices are even, there are seven types of
smoothings corresponding to $[G]$ (and seven types of smoothings
corresponding to $[G']$): in each of the three vertex we have two
possible smoothings, and one case is ruled out because of a free
loop (which yields at least two unicursal components, thus having no
impact in $[G]$ or $[G']$). When considering $G$, three of these
seven cases coincide (this triple is denoted by $1$), so, in $\ZG$
it remains exactly one of these two cases. Amongst the smoothings of
the diagram $G'$, the other three cases coincide (they are marked by
$2$). Thus, both in $[G]$ and $[G']$ there are five types of
summands marked by $1,2,3,4,5$.

These five cases are in one-to-one correspondence (see Fig.
\ref{razved31}) and they yield the equality  $[G]=[G']$).

\begin{figure}
\centering\includegraphics[width=200pt]{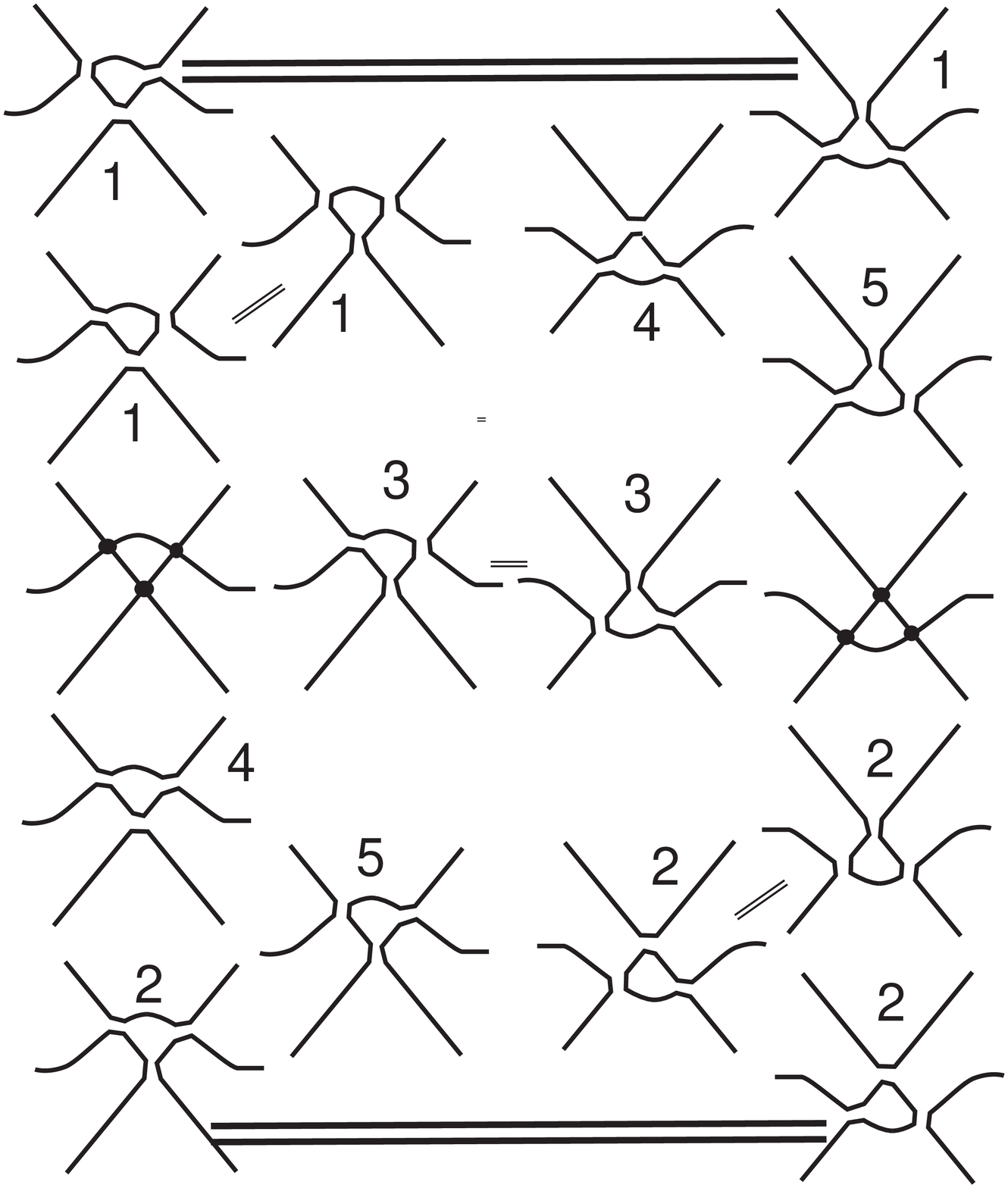}
\caption{Correspondences of smoothings with respect to $\Omega_3$
with three even vertices} \label{razved31}
\end{figure}

If amongst the three vertices taking part in $\Omega_3$ we  have
exactly one even vertices (say $a\to a'$), we get the situation
depicted in Fig. \ref{razved32}.

\begin{figure}
\centering\includegraphics[width=200pt]{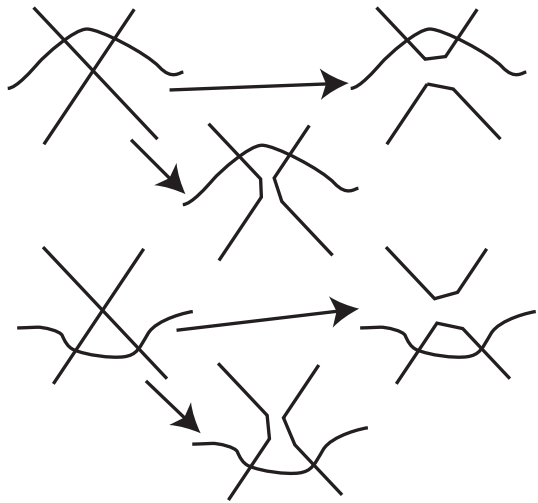}
\caption{Correspondence between smoothings for $\Omega_3$ with one
even vertex} \label{razved32}
\end{figure}

From this figure we see that those smoothings where $a$ (resp.,
$a'$) is smoothed {\em vertically}, give identical summands in $[G]$
and in $[G']$, and those smoothings where $a$ and $a'$ are smoothed
{\em horizontally}, are in one-to-one correspondence for $G$ and
$G'$, and the corresponding summands are obtained by applying two
second Reidemeister moves. This proves that $[G]=[G']$ in $\ZG$.

\end{proof}

\section{A functorial mapping}

It turns out that the odd axioms listed above lead to a simple and
powerful map on the set of free knots and, more generally, on the
set of virtual knots.

Let $K$ be a virtual knot diagram. Let $f$ be a diagram obtained
from $K$ by making all {\em odd} crossings virtual. In other words,
we remove all odd chords.

The following theorem follows from definitions.

\begin{thm}
$f$ is a well-defined map on the set of all virtual knots. For a
virtual knot diagram $K$, $f(K)=K$ iff the atom corresponding to $K$
is orientable. Otherwise, $f(K)$ has strictly less classical
crossings than $K$.
\end{thm}

This theorem leads to a new proof of a partial case the following
\begin{thm}[First proved by O.Ya. Viro and V.O.Manturov, 2005, first
published in \cite{IM}] The set of virtual links with orientable
atoms is closed. In other words, if two virtual diagrams $K$ and
$K'$ have orientable atoms and they are equivalent, then there is a
sequence of diagrams $K=K_{0}\to K_{1}\dots\to K_{n}=K'$ all having
orientable atoms where $K_{i}$ is obtained from $K_{i}$ by a
Reidemeister move.
\end{thm}

The case of this theorem for knots (1-component links) is then
proved as follows: take any sequence $K=L_{0}\to\dots\to L_{m}=K'$
and apply the map $f$ to it as many times as necessary (that is, to
make all intermediate knot diagrams having orientable atom). This
leads to the desired sequence.

Thus, odd chords of Gauss diagrams are holders of the ``atom
non-orientability'' condition for virtual knots (one-component
links). It would be interesting to find an appropriate ``parity''
for the case of arbitrary component links.

The map $f$ is also important because it gives rise to a natural
filtration on the set of virtual knots: ${\cal K_{0}}\subset {\cal
K_{1}}\in {\cal K_{2}}\subset \dots \subset {\cal K_{n}}$ where
${\cal K_{n}}$ denotes the set of all virtual knots $K$ where
$f^{n}(K)$ has a diagram with orientable atom. In particular, all
classical knots are in ${\cal K_{0}}$.

\section{Further applications}

Theorem \ref{mainthm} can be treated independently from knots. It
allows to establish non-equivalence of some quandratic forms over
$\Z_{2}$, which are considered up to a certain quite natural and
simple equivalence: (de)stabilizations, i.\,e., addition/removal of
a vector perpendicular to the whole space, addition/removal of a
pair of ``similar'' vectors, and a simple transformation
corresponding to the Reidemeister 3 move.

We construct an invariant of such quandratic forms valued in formal
${\bf Z}_{2}$-linear combinations of much simpler equivalence
classes of quadratic forms. These classes are actually classified by
intersection matrices of odd bunches.

At the level of virtual knots, Corollary \ref{sld} immediately
yields minimality of some virtual knot diagrams.

Actually, the only thing we need from knots and chord diagrams are
intersection graphs \cite{CDL}. Thus, the results given above
generalize for the case of {\em graph-links} by Ilyutko-Manturov
\cite{IM} and looped graphs by Traldi and Zulli \cite{ZT} which are
equivalence classes of some graphs modulo formal Reidemeister moves
(which are defined with respect to the information coming from the
intersection matrix).

The approach given above proves the {\em non-realizability} of some
looped graphs and {\em graph-links}.

It suffices to take the graph $P_5$, representing the $1$-frame of
the pentagon pyramid. This graph is irreducibly odd, thus any
corresponding looped graph in sense of \cite{TZ} has no realizable
representative: if there were any realizable representative, one
would easily obtain a realizability of a chord diagram whose
intersection graph contains $P_5$ as a subgraph.

The same argument works for graph-links \cite{IM} with {\em
non-orientable atoms}, \cite{MyBook}.

\end{document}